\documentclass{amsart}
\usepackage{amssymb,amscd,amsmath,latexsym,amsthm}
\usepackage[mathcal]{euscript}
\usepackage[small,nohug,heads=LaTeX]{diagrams}
\usepackage{float}
\diagramstyle[labelstyle=\scriptstyle]

\newtheorem{thm}{Theorem}[section]
\newtheorem*{thm1}{Theorem 1 (\ref{main.theorem})}
\newtheorem*{thm2}{Theorem 2 (\ref{general.quotients.iso.to.sl.lie.algebra} and \ref{SL.quotients})}
\newtheorem*{cor3}{Corollary 3 (\ref{homology.level.p.general})}
\newtheorem*{thm4}{Theorem 4 (\ref{congruence.differential})}
\newtheorem*{thm5}{Theorem 5 (\ref{infinitely.generated.cohomology})}
\newtheorem*{thm6}{Theorem 6 (\ref{bound.gaussian.cohomology})}

\newtheorem{cor}[thm]{Corollary}
\newtheorem{lem}[thm]{Lemma}
\newtheorem{prop}[thm]{Proposition}
\newtheorem{defin}[thm]{Definition}

\newcommand{\Z}{{\mathbb Z}}
\newcommand{\F}{{\mathbb F}}

\newcommand{\ilim}{\mathop{\varprojlim}\limits}

\DeclareMathOperator{\gr}{gr}
\DeclareMathOperator{\Mat}{Mat}

\numberwithin{equation}{section}

\begin{document}

\title{Lie Algebras and Cohomology of Congruence Subgroups for $SL_n(R)$}
\author{Jonathan Lopez}
\address{Department of Mathematics,
Niagara University, Niagara University, NY 14109}
\email{jlopez@niagara.edu}

\begin{abstract}
Let $R$ be a commutative ring that is free of rank $k$ as an abelian group, $p$ a prime, and $SL_n(R)$ the special linear group.  We show that the Lie algebra associated to the filtration of $SL_n(R)$ by $p$-congruence subgroups is isomorphic to the tensor product $\mathfrak{sl}_n(R\otimes_{\Z}\Z/p)\otimes_{\F_p}t\F_p[t]$, the Lie algebra of polynomials with zero constant term and coefficients $n\times n$ traceless matrices with entries polynomials in $k$ variables over $\F_p$.

We use the Lie algebra structure along with the Lyndon-Hochschild-Serre spectral sequence to compute the $d^2$ homology differential for certain central extensions involving quotients of $p$-congruence subgroups.  We also use the underlying group structure to obtain several homological results.  For example, we compute the first homology group of the level $p$-congruence subgroup for $n\geq3$.  We show that the cohomology groups of the level $p^r$-congruence subgroup are not finitely generated for $n=2$ and $R=\Z[t]$.  Finally, we show that for $n=2$ and $R=\Z[i]$ (the Gaussian integers) the second cohomology group of the level $p^r$-congruence subgroup has dimension at least two as an $\F_p$-vector space.
\end{abstract}


\maketitle

\section{Introduction}\label{introduction}
The method of constructing a Lie algebra from a filtered group is classical (see \cite{serre}, for example), and congruence subgroups have been widely studied for many years by mathematicians such as Bass, Milnor, and Serre (in \cite{BMS}), among others.  Linear groups over commutative rings admit natural filtrations by $p$-congruence subgroups, i.e., subgroups consisting of matrices that reduce to the identity modulo powers of $p$.  As a result, we are able to construct an associated Lie algebra, and the structure of this Lie algebra encodes homological information, as well as information concerning the structure of certain central extensions involving quotients of congruence subgroups.  It should be noted that filtrations of this type, i.e., filtrations of linear groups by $p$-congruence subgroups, appear in linearity problems of groups (see \cite{cohen} and \cite{lubotzky}).

We start with a commutative ring $R$ that is free of rank $k$ as a $\Z$-module (unless otherwise stated, we allow $k$ to be infinite).  Let $G_n(R)$ denote an $n\times n$ matrix group with coefficients in $R$.  In this paper, we assume $G_n(R)=SL_n(R)$ unless otherwise stated.  Examples involving other linear groups appear in \cite{lopez}.  For a prime $p$, there is a filtration 
\begin{align}\label{filtration}
\cdots\subseteq\Gamma(G_n(R),p^r)\subseteq\cdots\subseteq\Gamma(G_n(R),p^2)\subseteq\Gamma(G_n(R),p)
\end{align}
for $r\geq1$, where
\begin{align}
\begin{CD}
\Gamma(G_n(R),p^r)=\ker(G_n(R)@>>>G_n(R\otimes_{\Z}\Z/p^r)).
\end{CD}
\end{align}
Classically, there is a Lie algebra
\begin{align}
\gr_*(\Gamma(G_n(R),p))=\bigoplus_{r\geq1}{\Gamma(G_n(R),p^r)/\Gamma(G_n(R),p^{r+1})}
\end{align}
associated to filtration (\ref{filtration}) in which the Lie bracket
\begin{align}
\begin{CD}
[-,-]:\gr_*(\Gamma(G_n(R),p))\otimes_{\F_p}\gr_*(\Gamma(G_n(R),p))@>>>\gr_*(\Gamma(G_n(R),p))
\end{CD}
\end{align}
is induced by the restriction of the group commutator map in $G_n(R)$.  Our main results include an explicit description of the Lie algebra structure of $\gr_*(\Gamma(G_n(R),p))$ for the case $G_n(R)=SL_n(R)$, as well as the following theorem.

\begin{thm1}\label{main.theorem.introduction}
Let $\mathfrak{g}=\mathfrak{sl}_n(R\otimes_{\Z}\Z/p)$, and let $\mathfrak{L}$ be the kernel of the evaluation map $\begin{CD}\mathfrak{g}\otimes_{\F_p}\F_p[t]@>{t=0}>>\mathfrak{g}\end{CD}$.  Then
\begin{align}
\mathfrak{L}\cong\gr_*(\Gamma(G_n(R),p))
\end{align}
as Lie algebras.
\end{thm1}

We also explore the group structure of the filtration quotients for the filtration of $\Gamma(G_n(R),p)$ given in (\ref{filtration}).  The group structure is used to compute the structure of the Lie algebra $\gr_*(\Gamma(G_n(R),p))$, and it is also used in the proof of Theorem \ref{main.theorem.introduction}.  Our main results concerning the filtration quotients are as follows, where $V=\{v_i\}_{i\in I}$ is a $\Z$-basis for $R$ (recall that $R$ is a free $\Z$-module of rank $k$). 

\begin{thm2}\label{SL.quotients.main}
For $n\geq2$ and $r\geq1$,
\begin{align}
\Gamma(G_n(R),p^r)/\Gamma(G_n(R),p^{r+1})\cong\bigoplus_{n^2-1}{\F_p[V]}.
\end{align}
Moreover,
\begin{align}
\Gamma(G_n(R),p^r)/\Gamma(G_n(R),p^{r+1})\cong\mathfrak{sl}_n(\F_p[V])
\end{align}
as abelian groups.
\end{thm2}

The next corollary is a direct application of this theorem and a result of Bass-Milnor-Serre \cite{BMS}.

\begin{cor3}\label{homology}
For $n\geq3$,
\begin{align}
H_1(\Gamma(G_n(R),p);\Z)\cong\mathfrak{sl}_n(\F_p[V]).
\end{align}
\end{cor3}

Finally, we use the structure of $\gr_*(\Gamma(G_n(R),p))$ to obtain some homological information concerning $p$-congruence subgroups for $SL_n(R)$.  The Lie algebra structure can be used to explicitly compute the $d^2$ homology differential in the Lyndon-Hochschild-Serre spectral sequence for certain central extensions involving quotients of $p$-congruence subgroups.  More specifically, suppose one has a central extension
\begin{align}\label{differential.extension.intro}
\begin{CD}
1@>>>\Gamma_{r+s-1}/\Gamma_{r+s}@>>>\Gamma_r/\Gamma_{r+s}@>>>\Gamma_r/\Gamma_{r+s-1}@>>>1,
\end{CD}
\end{align}
where we have used the convention $\Gamma_r=\Gamma(G_n(R),p^r)$.  Furthermore, suppose that defining a map
\begin{align}\label{map.of.free.group}
\begin{CD}
\theta:F[x,y]@>>>\Gamma_r/\Gamma_{r+s}
\end{CD}
\end{align}
induces a map of extensions as in Figure \ref{induced.map.of.extensions.intro},
\begin{figure}[h!]
$$\begin{CD}
  1@>>>K@>>>F[x,y]@>>>\Z\oplus\Z @>>>1\\
  @. @VV{\theta}V @VV{\theta}V @VV{\theta}V @.\\
  1@>>>\Gamma_{r+s-1}/\Gamma_{r+s}@>>>\Gamma_r/\Gamma_{r+s}@>>>\Gamma_r/\Gamma_{r+s-1}@>>>1\\
\end{CD}$$
\caption{Map of Extensions Induced by (\ref{map.of.free.group})}\label{induced.map.of.extensions.intro}
\end{figure}
where $F[x,y]$ is the free group on two generators.  In many cases, the Lie algebra structure can be used to explicitly compute $d^2$ for extension (\ref{differential.extension.intro}), thus informing on $H_1(\Gamma_r/\Gamma_{r+s})$.  The main result concerning this is presented in the next theorem.

\begin{thm4}\label{differential.thm}
Let $\chi$ be the generator of $H_2(\Z\oplus\Z)\cong\Z$.  Given the map of extensions $\theta$ described above,
\begin{align}
d^2(\theta_*(\chi))=[\theta_*(x),\theta_*(y)].
\end{align}
\end{thm4}

We are also able to obtain some cohomological information by examining certain subgroups of $\Gamma(G_n(R),p^r)$.  For the case $n=2$ and $R=\Z[t]$, we are able to produce a subgroup $\bigoplus_{i\geq0}\Z\subseteq\Gamma(G_2(R),p^r)$, from which we are able to prove the following result.

\begin{thm5}
Suppose $n=2$ and $R=\Z[t]$.  For all $j,r\geq1$ and for all primes $p$,
\begin{align}
H^j(\Gamma(G_2(\Z[t]),p^r);\F_p)
\end{align}
is infinitely generated.
\end{thm5}

For the case $n=2$ and $R=\Z[i]$, we construct a map
\begin{align}
\begin{CD}
(\Z\oplus\Z)*(\Z\oplus\Z)@>>>\Gamma(G_2(R),p^r)
\end{CD}
\end{align}
that can be used to show that $H^2(\Gamma(G_2(R),p^r);\F_p)$ has dimension at least 2 as an $\F_p$-vector space.  In particular, this confirms that the level $p^r$-congruence subgroup for $SL_2(\Z[i])$ is not free for any $r\geq1$.

\begin{thm6}
Suppose $n=2$ and $R=\Z[i]$.  For all $r\geq1$ and for all primes $p$,
\begin{align}
H^2(\Gamma(G_2(\Z[i]),p^r);\F_p)\supseteq\F_p\oplus\F_p.
\end{align}
\end{thm6}

\noindent\textit{Acknowledgments.}  The author would like to thank Fred Cohen for suggesting many of the problems discussed in this paper and for countless conversations concerning these topics over the past several years.  The author would also like to thank Stratos Prassidis for helpful suggestions concerning the layout of this paper, as well as several insightful discussions concerning its content.

\section{Preliminaries, Basic Notation}

We now collect the notation that will be used throughout this paper: $R$ is a commutative ring that is free of rank $k$ as an abelian group; $V=\{v_i\}_{i\in I}$ is a $\Z$-basis for $R$; $G_n(R)$ is an $n\times n$ matrix group with coefficients in $R$; $p$ is a prime; and $\Gamma(G_n(R),p^r)$ is the kernel of the mod-$p^r$ reduction map
\begin{align}
\begin{CD}
G_n(R)@>>>G_n(R\otimes_{\Z}\Z/p^r).
\end{CD}
\end{align}
The following is a classical result (see, for example, \cite{dusautoy}):

\begin{thm}\label{lie.algebra.thm}
There is a Lie algebra $\gr_{*}(\Gamma(G_n(R),p))$ associated to the filtration 
\begin{align}\label{filtration.thm}
\cdots\subseteq\Gamma(G_n(R),p^r)\subseteq\cdots\subseteq\Gamma(G_n(R),p^2)\subseteq\Gamma(G_n(R),p)
\end{align}
in which the Lie bracket is induced by the restriction of the commutator map in $G_n(R)$,
\begin{align}
\begin{CD}
[-,-]:\Gamma(G_n(R),p^r)\times\Gamma(G_n(R),p^s)@>>>\Gamma(G_n(R),p^{r+s}).
\end{CD}
\end{align}
\end{thm} 
As one of the main results of this paper will be describing the construction and structure of this Lie algebra for $SL_n(R)$, it may be worthwhile to first consider an elementary case ($n=2$, $R=\Z$) in order to explicitly see the full calculation of the Lie algebra structure.

\begin{thm}\label{SL(2,Z)}
Suppose $n=2$ and $R=\Z$, so that $G_n(R)=SL_2(\Z)$.  For $r\geq1$, consider the following matrices:
\begin{align}
A_{12,r}=\begin{pmatrix}1&p^r\\0&1\\\end{pmatrix}\qquad A_{21,r}=\begin{pmatrix}1&0\\p^r&1\\\end{pmatrix}\qquad A_{11,r}=\begin{pmatrix}1+p^r&0\\0&1-p^r\\\end{pmatrix}.
\end{align}
The Lie algebra $\gr_{*}(\Gamma(SL_2(\Z),p))$ is generated (as a restricted Lie algebra over $\F_p$) by the three matrices $A_{12,1}$, $A_{21,1}$, and $A_{11,1}$.  Furthermore, the following relations are satisfied for all $r,s\geq1$ and for all primes $p$:
\begin{enumerate}
\item $[A_{11,r},A_{12,s}]=A_{12,r+s}^2$
\item $[A_{11,r},A_{21,s}]=A_{21,r+s}^{-2}$
\item $[A_{12,r},A_{21,s}]=A_{11,r+s}$
\item $A_{ij,r}^p=A_{ij,r+1}$ for $(i,j)\in\{(1,2),(2,1),(1,1)\}$.
\end{enumerate}
\end{thm}

One feature of the Lie algebra $\gr_{*}(\Gamma(G_n(R),p))$ is the fact that it is finitely generated as a restricted Lie algebra whenever $k<\infty$.  This is suggested by Theorem \ref{SL(2,Z)} for the case $n=2$ and $k=1$, and is, in fact, true for all $n\geq2$ and $k<\infty$.

We now provide a description of the Lie algebra $\mathfrak{L}$ introduced in Theorem \ref{main.theorem.introduction}.  Recall that $\mathfrak{sl}_n(\F_p[V])$ is the Lie algebra of $n\times n$ traceless matrices over $\F_p[V]$, where the Lie bracket is defined by
\begin{align}
[A,B]=AB-BA
\end{align}
for $A,B\in\mathfrak{sl}_n(\F_p[V])$.  Let $\F_p[t]$ be the polynomial ring in one indeterminate over $\F_p$, and let $I$ denote the kernel of the evaluation map
\begin{align}\label{eval.map}
\begin{CD}
  \F_p[t]@>{t=0}>>\F_p.
\end{CD}
\end{align}
Notice that $I$ is the $\F_p$-linear span of $\{t^i\,|\,i>0\}$.

Writing $\mathfrak{g}=\mathfrak{sl}_n(\F_p[V])$, consider the Lie algebra $\mathfrak{g}\otimes_{\F_p}\F_p[t]$, where the Lie bracket is defined by
\begin{align}\label{lie.bracket.sl.tensor.poly}
[A\otimes t^i,B\otimes t^j]=[A,B]\otimes t^{i+j}=(AB-BA)\otimes t^{i+j}
\end{align}
for $A,B\in\mathfrak{g}$ and $i,j\geq0$.  The evaluation map (\ref{eval.map}) induces a morphism of Lie algebras 
\begin{align}
\begin{CD}
  \mathfrak{g}\otimes_{\F_p}\F_p[t]@>{t=0}>>\mathfrak{g}
\end{CD}
\end{align}
whose kernel is exactly $\mathfrak{L}=\mathfrak{g}\otimes_{\F_p}I$.  It is easy to verify that $\mathfrak{L}$ is a Lie algebra with Lie bracket as defined in \ref{lie.bracket.sl.tensor.poly}.  These facts are collected in the next lemma.

\begin{lem}
The split short exact sequence of abelian groups
\begin{align}
\begin{CD}
0@>>>I@>>>\F_p[t]@>{t=0}>>\F_p@>>>0
\end{CD}
\end{align}
gives rise to a split short exact sequence of Lie algebras
\begin{align}\label{ses.lie.algebras}
\begin{CD}
0@>>>\mathfrak{L}@>>>\mathfrak{g}\otimes_{\F_p}\F_p[t]@>{t=0}>>\mathfrak{g}@>>>0.
\end{CD}
\end{align}
\end{lem}

\begin{proof}
The splitting $\begin{CD}\F_p@>>>\F_p[t]\end{CD}$ in the first short exact sequence is given by the inclusion map $x\mapsto x$.

Since the tensor product is over the field $\F_p$, the fact that (\ref{ses.lie.algebras}) is exact is immediate.  The map $\begin{CD}\mathfrak{g}@>>>\mathfrak{g}\otimes_{\F_p}\F_p[t]\end{CD}$ defined by $x\mapsto x\otimes t^0$ gives the splitting.  In fact, since $[x,y]\mapsto[x,y]\otimes t^0=[x\otimes t^0,y\otimes t^0]$, the sequence is split as Lie algebras.
\end{proof}

\begin{thm}\label{main.theorem}
Let $\mathfrak{g}=\mathfrak{sl}_n(\F_p[V])$, and let $\mathfrak{L}$ be the kernel of the evaluation map
\begin{align}
\begin{CD}
\mathfrak{g}\otimes_{\F_p}\F_p[t]@>{t=0}>>\mathfrak{g}
\end{CD}
\end{align}
as discussed above.  Then
\begin{align}
\mathfrak{L}\cong\gr_*(\Gamma(G_n(R),p))
\end{align}
as Lie algebras.
\end{thm}

An immediate application of Theorem \ref{main.theorem} is the computation of the Lie algebra homology of $\gr_*(\Gamma(G_n(R),p))$.  Let $\left\langle t^i\right\rangle$ denote the $\F_p$-span of $t^i$.

\begin{cor}
Let $\mathfrak{g}=\mathfrak{sl}_n(\F_p[V])$ as above.  Then
\begin{align}
H_1(\gr_*(\Gamma(G_n(R),p)))=(\mathfrak{g}\otimes_{\F_p}\left\langle t\right\rangle)\oplus(\bigoplus_{i\geq2}{H_1(\mathfrak{g})\otimes_{\F_p}\left\langle t^i\right\rangle}).
\end{align}
\end{cor}

The proof of Theorem \ref{main.theorem} is included in Section \ref{proof.of.main.theorem}.  We first record some basic results concerning the group structure of the filtration quotients for the filtration of $\Gamma(G_n(R),p)$ defined in (\ref{filtration}) and (\ref{filtration.thm}).

\section{The structure of the filtration quotients for $\Gamma(G_n(R),p)$}

The following theorem is due to Lee and Szczarba \cite{lee.szczarba}.

\begin{thm}\label{homology.level.p}
Let $G_n(R)=SL_n(\Z)$, and suppose $n\geq3$.  There is an epimorphism
\begin{align}
\begin{CD}
\varphi_1:\Gamma(SL_n(\Z),p)@>>>\mathfrak{sl}_n(\F_p)
\end{CD}
\end{align}
whose kernel is the commutator subgroup.  Thus, $H_1(\Gamma(SL_n(\Z),p);\Z)\cong\mathfrak{sl}_n(\F_p)$.
\end{thm}

Generalizing to the level $p^r$-congruence subgroup and to coefficient rings other than $\Z$, we obtain the following result.

\begin{thm}\label{general.quotients.iso.to.sl.lie.algebra}
Suppose $n\geq2$ and $r\geq1$.  For $G_n(R)=SL_n(R)$, there is an epimorphism
\begin{align}
\begin{CD}
\varphi_r:\Gamma(G_n(R),p^r)@>>>\mathfrak{sl}_n(\F_p[V])
\end{CD}
\end{align}
whose kernel is $\Gamma(G_n(R),p^{r+1})$.  Thus,
\begin{align}\Gamma(G_n(R),p^r)/\Gamma(G_n(R),p^{r+1})\cong\mathfrak{sl}_n(\F_p[V]).
\end{align}
\end{thm}
The proof of this theorem will be included in Section \ref{general.quotients.iso.to.sl.lie.algebra.proof}.  An immediate application of Theorem \ref{general.quotients.iso.to.sl.lie.algebra} is the following corollary, which uses a result of Bass-Milnor-Serre \cite{BMS}.

\begin{cor}\label{homology.level.p.general}
For $n\geq3$ and $G_n(R)=SL_n(R)$,
\begin{align}
\Gamma(G_n(R),p^2)=[\Gamma(G_n(R),p),\Gamma(G_n(R),p)].
\end{align}
In particular,
\begin{align}
H_1(\Gamma(G_n(R),p);\Z)\cong\mathfrak{sl}_n(\F_p[V]).
\end{align}
\end{cor}

\begin{proof}
That $\Gamma(G_n(R),p^2)\supseteq[\Gamma(G_n(R),p),\Gamma(G_n(R),p)]$ follows from the fact that filtration (\ref{filtration.thm}) is Lie-like, as will be proved in Lemma \ref{SL.commutator}.

Let $e_{ij}$ denote the $n\times n$ matrix with a 1 in the $(i,j)$ position and zeros elsewhere.  Let $E(G_n(R),p^r)$ be the normal subgroup of $\Gamma(G_n(R),p^r)$ generated by matrices of the form $1+\alpha e_{ij}$, where $i\neq j$, 1 denotes the identity matrix, $\alpha\in R$, and $\alpha\equiv0\mod{p^r}$.  By the work of Bass-Milnor-Serre \cite{BMS}, we have $\Gamma(G_n(R),p^r)=E(G_n(R),p^r)$ and
\begin{align}
E(G_n(R),p^{r+s})\subseteq[E(G_n(R),p^r),E(G_n(R),p^s)]
\end{align}
for $n\geq3$.  Thus,
\begin{align}
\Gamma(G_n(R),p^2)\subseteq[E(G_n(R),p),E(G_n(R),p)]\subseteq[\Gamma(G_n(R),p),\Gamma(G_n(R),p)].
\end{align}
It  follows that $\Gamma(G_n(R),p^2)=[\Gamma(G_n(R),p),\Gamma(G_n(R),p)]$ for $n\geq3$.  Using this along with Theorem \ref{general.quotients.iso.to.sl.lie.algebra} (for $r=1$), $H_1(\Gamma(G_n(R),p);\Z)\cong\mathfrak{sl}_n(\F_p[V])$ for $n\geq3$.
\end{proof}

Next, we provide an alternate description of the filtration quotients that we use to describe the structure of the Lie algebra $\gr_*(\Gamma(G_n(R),p))$.  We make use of the following fact, whose proof we include here for completeness.

\begin{lem}\label{determinant.1+p^rA}
Let $p$ be a prime and suppose $A=(a_{ij})\in\Mat_n(R)$.  For $r\geq1$,
\begin{align}
\det(1+p^rA)\equiv(1+p^rtr(A))\mod p^{r+1},
\end{align}
where $tr(A)$ denotes the trace of $A$.  In particular, if $\det(1+p^rA)\equiv1\mod p^{r+1}$, it must be the case that
\begin{align}
tr(A)\equiv0\mod p.
\end{align}
In other words, $a_{nn}\equiv-(a_{11}+\cdots+a_{n-1,n-1})\mod p$.
\end{lem}

\begin{proof} 
The proof is by induction on $n$.  For $n=2$, suppose 
\begin{align}
1+p^rA=\begin{pmatrix}1+p^ra_{11}&p^ra_{12}\\p^ra_{21}&1+p^ra_{22}\end{pmatrix},
\end{align}
where $a_{ij}\in R$.  Then 
\begin{align*}
\det(1+p^rA) & =       1+p^r(a_{11}+a_{22})+p^{2r}(a_{11}a_{22}-a_{12}a_{21})\\[.1em]
             & \equiv (1+p^r(a_{11}+a_{22}))\mod p^{r+1}.
\end{align*}
If $\det(1+p^rA)\equiv1\mod{p^{r+1}}$, it must be the case that $tr(A)\equiv0\mod p$.  This proves the claim for the case $n=2$.

Suppose $n-1$ and smaller cases have been proved.  Let $\hat{A}_{ij}$ be the $(n-1)\times(n-1)$ matrix obtained by removing the $i^{\text{th}}$ row and $j^{\text{th}}$ column from $1+p^rA=1+p^r(a_{ij})$.  By induction, we can write
\begin{align*}
\det(1+p^rA)&=(1+p^ra_{11})\det(\hat{A}_{11})+p^r\sum_{j=2}^n{(-1)^{j+1}a_{1j}\det(\hat{A}_{1j})}\\[.1em]
&\equiv[(1+p^ra_{11})(1+p^r(a_{22}+\cdots+a_{nn}))\\[.1em]
&\qquad\qquad+p^r\sum_{j=2}^n{(-1)^{j+1}a_{1j}\det(\hat{A}_{1j})}]\mod p^{r+1}.
\end{align*}
Notice that each of $\hat{A}_{1j}$ for $2\leq j\leq n$ is an $(n-1)\times(n-1)$ matrix whose first column is
\begin{align}
\begin{pmatrix}p^ra_{21}\\p^ra_{31}\\\vdots\\p^ra_{nn}\end{pmatrix}.
\end{align}
If we use the cofactor expansion on this column to compute the determinant of $\hat{A}_{1j}$, it is clear that $\det(\hat{A}_{1j})\equiv0\mod p^r$ for $2\leq j\leq n$, say $\det(\hat{A}_{1j})=p^r\alpha_{1j}$.  Using this above, we see that 
\begin{align*}
\det(1+p^rA)&\equiv[(1+p^ra_{11})(1+p^r(a_{22}+\cdots+a_{nn}))\\[.1em]
&\qquad\qquad+p^{2r}\sum_{j=2}^n{(-1)^{j+1}a_{1j}\alpha_{1j}}]\mod p^{r+1}\\[.1em]
&\equiv[(1+p^ra_{11})(1+p^r(a_{22}+\cdots+a_{nn}))]\mod p^{r+1}\\[.1em]
&\equiv[1+p^rtr(A)+p^{2r}a_{11}(a_{22}+\cdots+a_{nn})]\mod p^{r+1}\\[.1em]
&\equiv[1+p^rtr(A)]\mod p^{r+1}.
\end{align*}
If $\det(1+p^rA)\equiv1\mod{p^{r+1}}$, it must be the case that $tr(A)\equiv0\mod p$.
This completes the proof of the lemma.
\end{proof}

\begin{thm}\label{SL.s.stage.quotients}
Let $V$ be a $\Z$-basis for $R$.  Suppose $|V|<\infty$ or $s=1$.  For $r\geq s\geq1$ and $p$ a prime,
\begin{align}
\Gamma(G_n(R),p^r)/\Gamma(G_n(R),p^{r+s})\cong\bigoplus_{n^2-1}{\Z/p^s\Z[V]}.
\end{align}
\end{thm}
The proof is by induction on $s$.  For the sake of continuity, we prove the case $s=1$ in the next lemma, and complete the induction in Section \ref{SL.s.stage.quotients.proof}.  We first state a corollary that follows directly from Theorem \ref{SL.s.stage.quotients}.

\begin{cor}
For all primes $p$,
\begin{align}
\ilim_s{\Gamma(G_n(R),p^s)/\Gamma(G_n(R),p^{2s})}=\bigoplus_{n^2-1}\hat{\Z}_p[V],
\end{align}
where $\hat{\Z}_p$ denotes the $p$-adic integers.
\end{cor}

\begin{proof} Follows immediately from Theorem \ref{SL.s.stage.quotients} and the fact that inverse limits commute with direct sums.
\end{proof}

\begin{lem}\label{SL.quotients}
For $r\geq1$ and $p$ a prime,
\begin{align}
\Gamma(G_n(R),p^r)/\Gamma(G_n(R),p^{r+1})\cong\bigoplus_{n^2-1}{\F_p[V]},
\end{align}
where $V$ is a $\Z$-basis for $R$.
\end{lem}

\begin{proof} Consider the commutative diagram in Figure \ref{comm.diag},
\footnotesize
\begin{figure}[h!]
$$\begin{CD}
  &&1&&1&&1&&\\
  && @VVV @VVV @VVV && \\
  1@>>>\Gamma(G_n(R),p^{r+1})@>{Id}>>\Gamma(G_n(R),p^{r+1})@>>>\{1\}@>>>1 \\
  && @VVV @VVV @VVV && \\
  1@>>>\Gamma(G_n(R),p^r)@>>>SL_n(R)@>{\phi_r}>>SL_n(R\otimes_{\Z}\Z/p^r)@>>>1 \\
  && @VVV @VV{\phi_{r+1}}V @VV{Id}V && \\
  1@>>>\text{ker}\,\theta_r@>>>SL_n(R\otimes_{\Z}\Z/p^{r+1})@>{\theta_r}>>SL_n(R\otimes_{\Z}\Z/p^r)@>>>1 \\
  && @VVV @VVV @VVV && \\
  &&1&&1&&1&&
\end{CD}$$
\caption{}\label{comm.diag}
\end{figure}
\normalsize
where $\phi_r$, $\phi_{r+1}$, and $\theta_r$ are the natural reduction maps and $\phi_r=\theta_r\circ\phi_{r+1}$.  In \cite{shimura}, it is shown that reduction maps such as $\phi_r$, $\phi_{r+1}$, and $\theta_r$ are surjections.  Given this, it is immediate that the three rows and two right columns are exact.  A standard diagram chase can be used to show that the first column is also exact.  This gives
\begin{align}
\text{ker}\,\theta_r\cong\Gamma(G_n(R),p^r)/\Gamma(G_n(R),p^{r+1}).
\end{align}
The reader should note that this isomorphism allows us to consider
\begin{align}
\Gamma(G_n(R),p^r)/\Gamma(G_n(R),p^{r+1})\subseteq SL_n(R\otimes_{\Z}\Z/p^{r+1}).
\end{align}
This will prove to be a useful fact at many points throughout this paper.

A priori, an element of $\text{ker}\,\theta_r$ is a matrix of the form $1+p^rA$ where $A=(a_{ij})\in\Mat_n(R)$ and $\det(1+p^rA)\equiv1\mod p^{r+1}$.  Since $1+p^rA\in SL_n(R\otimes_{\Z}\Z/p^{r+1})$, we can assume $A\in\Mat_n(R\otimes_{\Z}\Z/p)$.

Define a map $\begin{CD}\Phi:\bigoplus_{n^2-1}{\F_p[V]}@>>>\text{ker}\,\theta_r\end{CD}$ by
\begin{align}
(a_{11},\ldots,a_{1n},a_{21},\ldots,a_{2n},\ldots,a_{n1},\ldots,a_{n,n-1})\mapsto1+p^r(a_{ij}),
\end{align}
where $a_{nn}=-(a_{11}+a_{22}+\cdots+a_{n-1,n-1})$.  By Lemma \ref{determinant.1+p^rA}, we see immediately that $\Phi$ is surjective.  By inspection, the kernel is trivial, so that $\Phi$ is a bijection.  

To check that $\Phi$ is a homomorphism, notice that 
\begin{align*}
(1+p^rA)(1+p^rB) & =       1+p^r(A+B)+p^{2r}AB\\[.1em]
                 & \equiv  1+p^r(A+B)\mod p^{r+1}.
\end{align*}
Thus, $\Phi$ is the required isomorphism.  This completes the proof.
\end{proof}

This lemma has several immediate consequences.  Recall that $V=\{v_i\}_{i\in I}$ is a $\Z$-basis for $R$.  For $r\geq1$, $1\leq i,j\leq n$, and $i+j<2n$, define a family of matrices by
\begin{align}\label{generators}
A_{ij,k,r}=\begin{cases} 1+p^rv_ke_{ij}          &\text{if $i\neq j$}\\
                           1+p^rv_k(e_{ii}-e_{nn}) &\text{if $i=j$},\\
             \end{cases}
\end{align}
where 1 denotes the $n\times n$ identity matrix, and the matrices $\{e_{ij}\}$ are as defined in the proof of Corollary \ref{homology.level.p.general}.

\begin{cor}\label{SL.generators}
For $r\geq1$, $\Gamma(G_n(R),p^r)/\Gamma(G_n(R),p^{r+1})$ is generated by the family of matrices $\{A_{ij,k,r}\}$ defined in \ref{generators}.
\end{cor}

\noindent\textit{Remark.}  It is clear that $\det(A_{ij,k,r})=1$ for $i\neq j$.  When $i=j$, notice that 
\begin{align}
\det(A_{ii,k,r})=(1+p^rv_k)(1-p^rv_k)=1-p^{2r}v_k^2\equiv1\mod p^{r+1}.
\end{align}
Since $\Gamma(G_n(R),p^r)/\Gamma(G_n(R),p^{r+1})\subseteq SL_n(R\otimes_{\Z}\Z/p^{r+1})$, this suffices.

\section{The Lie algebra associated to $SL_n(R)$}

We begin with some preliminary definitions concerning the existence of the associated Lie algebra, and then describe its structure completely.  The following definition is stated as in \cite{cohen}.

\begin{defin}
A filtration $\{F_n(G)\}_{n\geq0}$ of the group $G$ is said to be Lie-like if, for all $m,n\geq0$, the commutator map
\begin{align}
\begin{CD}
[-,-]:\,G\times G@>>>G
\end{CD}
\end{align}
given by $[g,h]\mapsto g^{-1}h^{-1}gh$ restricts to a map
\begin{align}\label{restriction.of.commutator}
\begin{CD}
[-,-]:\,F_m(G)\times F_n(G)@>>>F_{m+n}(G).
\end{CD}
\end{align}
\end{defin}

The next theorem is classical and can be found in \cite{serre}, for example.  The theorem gives a general procedure for constructing a Lie algebra from a filtered group whenever the filtration is Lie-like.

\begin{thm}\label{serre.thm}
Given a Lie-like filtration for the group $G$, there is an associated Lie algebra
\begin{align}
\gr_{*}(G)=\bigoplus_{n\geq1}{F_n(G)/F_{n+1}(G)}.
\end{align}
The Lie bracket is obtained by linearly extending the maps
\begin{align}
\begin{CD}
[-,-]:\,F_m(G)/F_{m+1}(G)\times F_n(G)/F_{n+1}(G)@>>>F_{m+n}(G)/F_{m+n+1}(G)
\end{CD}
\end{align}
that are induced by the restriction of the commutator map.
\end{thm}

\noindent\textit{Remark.}  Note that these Lie algebras are not, in general, graded Lie algebras.  In other words, it is not necessarily true that $[x,y]=(-1)^{|x||y|}[y,x]$, and the graded Jacobi identity need not be satisfied.

\vskip .1in

Next, we verify that the filtration of $\Gamma(G_n(R),p)$ given by (\ref{filtration.thm}) is Lie-like.

\begin{lem}\label{SL.commutator}
Suppose that $X\in\Gamma(G_n(R),p^r)$ and $Y\in\Gamma(G_n(R),p^s)$ for some $r,s\geq1$.  Then
$[X,Y]\in\Gamma(G_n(R),p^{r+s}),$
where $[X,Y]=X^{-1}Y^{-1}XY$ is the group commutator.
\end{lem}

\begin{proof}
We can write $X=1+p^rA$ and $Y=1+p^sB$ for some matrices $A,B\in\Mat_n(R)$.  Then
\begin{align*}
[1+p^rA,1+p^sB]&=(1+p^rA)^{-1}(1+p^sB)^{-1}(1+p^rA)(1+p^sB)\\[.1em]
&=\sum_{i=0}^{\infty}(-p^rA)^i\sum_{j=0}^{\infty}(-p^sB)^j(1+p^rA)(1+p^sB)\\[.1em]
&=1+p^{r+s}(AB-BA)+p^{r+s+1}C,
\end{align*}
where $C\in\Mat_n(R)$ is some combination of $A$ and $B$.  Thus,
\begin{align}
[X,Y]\in\Gamma(G_n(R),p^{r+s}),
\end{align}
as desired.
\end{proof}

By Theorem \ref{serre.thm}, there is an associated Lie algebra
\begin{align}
\gr_*(\Gamma(G_n(R),p))=\bigoplus_{r\geq1}{\Gamma(G_n(R),p^r)/\Gamma(G_n(R),p^{r+1})},
\end{align}
in which the Lie bracket is induced by the restriction of the commutator map.  In the next corollary, we use the notation $\Gamma_r=\Gamma(G_n(R),p^r)$.

\begin{cor}
The commutator map 
\begin{align}
\begin{CD}
[-,-]:\Gamma(G_n(R),p^r)\times\Gamma(G_n(R),p^s)@>>>\Gamma(G_n(R),p^{r+s})
\end{CD}
\end{align}
extends to a well-defined map on the filtration quotients
\begin{align}\label{filtration.quotient.bracket}
\begin{CD}
[-,-]:\Gamma_r/\Gamma_{r+1}\times\Gamma_s/\Gamma_{s+1}@>>>\Gamma_{r+s}/\Gamma_{r+s+1}.
\end{CD}
\end{align}
The Lie bracket
\begin{align}
\begin{CD}
[-,-]:\gr_*(\Gamma(G_n(R),p))\otimes_{\F_p}\gr_*(\Gamma(G_n(R),p))@>>>\gr_*(\Gamma(G_n(R),p))
\end{CD}
\end{align}
is obtained by extending the maps in (\ref{filtration.quotient.bracket}) linearly.
\end{cor}

\begin{proof}	
This follows immediately from the proof of Lemma \ref{SL.commutator} and \cite{serre}.
\end{proof}

The next result will be used to describe $\gr_*(\Gamma(G_n(R),p))$ as a restricted Lie algebra.  Recall that $\Gamma(G_n(R),p^r)/\Gamma(G_n(R),p^{r+1})$ for $r\geq1$ is generated by the family of matrices $\{A_{ij,k,r}\}$ defined in Corollary \ref{SL.generators}.

\begin{prop}\label{SL.p.power}
The Frobenius map
\begin{align}
\begin{CD}
\Gamma(G_n(R),p^r)@>>>\Gamma(G_n(R),p^{r+1})
\end{CD}
\end{align}
defined by 
$A_{ij,k,r}\mapsto(A_{ij,k,r})^p$ induces a map
\begin{align}
\begin{CD}
\psi_r^p:\Gamma(G_n(R),p^r)/\Gamma(G_n(R),p^{r+1})@>>>\Gamma(G_n(R),p^{r+1})/\Gamma(G_n(R),p^{r+2}).
\end{CD}
\end{align}
This map is an isomorphism, except in the case $p=2$, $r=1$.  Furthermore, $\psi_r^p(A_{ij,k,r})=A_{ij,k,r+1}$.
\end{prop}

\begin{proof}
Assume that $p\ne2$ or $r>1$.

A typical element of $\Gamma(G_n(R),p^r)/\Gamma(G_n(R),p^{r+1})$ is of the form $1+p^rA$ where $A\in\Mat_n(R)$.  In the proof of Lemma \ref{SL.quotients}, it was noted that
\begin{align}
\Gamma(G_n(R),p^r)/\Gamma(G_n(R),p^{r+1})\subseteq SL_n(R\otimes_{\Z}\Z/p^{r+1}),
\end{align}
so we can assume that $A\in\Mat_n(R\otimes_{\Z}\Z/p)$.  To see that the image of $\psi_r^p$ lies in $\Gamma(G_n(R),p^{r+1})/\Gamma(G_n(R),p^{r+2})$, notice that
\begin{align*}
\psi_r^p(1+p^rA) & =      (1+p^rA)^p\\[.1em]
                 & =      \sum_{i=0}^p{{\binom{p}{i}}(p^rA)^i}\\[.1em]
                 & \equiv (1+p^{r+1}A)\mod p^{r+2}.
\end{align*}

To check that $\psi_r^p$ is a homomorphism, consider two elements $1+p^rA$, $1+p^rB\in\Gamma(G_n(R),p^r)/\Gamma(G_n(R),p^{r+1})$.  Then
\begin{align*}
\psi_r^p((1+p^rA)(1+p^rB)) & =      \psi_r^p(1+p^r(A+B))\\[.1em]
                           & =      (1+p^r(A+B))^p\\[.1em]
                           & \equiv (1+p^{r+1}(A+B))\mod p^{r+2}\\[.1em]
                           & \equiv (1+p^{r+1}A)(1+p^{r+1}B)\mod p^{r+2}\\[.1em]
                           & =      \psi_r^p(1+p^rA)\psi_r^p(1+p^rB).
\end{align*}
Thus, $\psi_r^p$ is a homomorphism.

To see that $\psi_r^p$ is surjective, we show that it maps onto the generators of $\Gamma(G_n(R),p^{r+1})/\Gamma(G_n(R),p^{r+2})$.  For $i\neq j$,
\begin{align*}
(A_{ij,k,r})^p & =      (1+p^rv_ke_{ij})^p\\[.1em]
               & =      \sum_{l=0}^p{{\binom{p}{l}}(p^rv_ke_{ij})^l}\\[.1em]
               & \equiv (1+p^{r+1}v_ke_{ij})\mod p^{r+2}\\[.1em]
               & =      A_{ij,k,r+1}.
\end{align*}
For $i=j$, a similar computation gives the following:
\begin{align*}
(A_{ii,k,r})^p & =      (1+p^rv_k(e_{ii}-e_{nn}))^p \\[.1em]
               & =      \sum_{l=0}^p{{\binom{p}{l}}(p^rv_k(e_{ii}-e_{nn}))^l}\\[.1em]
               & \equiv (1+p^{r+1}v_k(e_{ii}-e_{nn}))\mod p^{r+2}\\[.1em]
               & =      A_{ii,k,r+1}.
\end{align*}
Thus, $\psi_r^p(A_{ij,k,r})=A_{ij,k,r+1}$, which gives that $\psi_r^p$ is surjective.

Finally, suppose that $\psi_r^p(1+p^rA)=1$.  Then it must be the case that
\begin{align}
p^{r+1}A\equiv0\mod p^{r+2}.
\end{align}
Since $A\in\Mat_n(R\otimes_{\Z}\Z/p)$, we conclude that $A\equiv0\mod p$.  Hence,
\begin{align}
(1+p^rA)\equiv1\mod p^{r+1},
\end{align}
so that $\psi_r^p$ is injective.  This completes the proof of the proposition.
\end{proof}

We are now ready to give the structure of the Lie algebra $\gr_*(\Gamma(G_n(R),p))$, and we do so in the next theorem.  In our computations, we make use of the following well-known result.

\begin{lem}\label{eij.bracket}
For the matrices $\{e_{ij}\}$ defined in the proof of Corollary \ref{homology.level.p.general}, we have the following:
\begin{align}
[e_{ij},e_{kl}]=e_{ij}e_{kl}-e_{kl}e_{ij}=\begin{cases}-e_{kj}        &\text{if $i=l$ and $j\neq k$,}\\
                                                        e_{il}        &\text{if $i\neq l$ and $j=k$,}\\
                                                        e_{il}-e_{kj} &\text{if $i=l$ and $j=k$,}\\
                                                        0             &\text{otherwise.}
                                          \end{cases}
\end{align}
\end{lem}

\begin{proof} It is clear that the first term, $e_{ij}e_{kl}$, is non-zero only when $j=k$, and the second term, $e_{kl}e_{ij}$, is non-zero precisely when $i=l$.  In these cases, $e_{ij}e_{kl}=e_{il}$ and $e_{kl}e_{ij}=e_{kj}$.  The result follows.
\end{proof}

In the next theorem, we use the notation $A_{ij,q_1q_2,r}=1+p^rv_{q_1}v_{q_2}e_{ij}$.  Here, $v_{q_1},v_{q_2}\in V$.

\begin{thm}\label{SL(n,Z[V]).lie.algebra}
Suppose $r\geq1$ and $p$ is a prime.  For the Lie algebra
\begin{align}
\gr_{*}(\Gamma(G_n(R),p))=\bigoplus_{r\geq1}{\Gamma(G_n(R),p^r)/\Gamma(G_n(R),p^{r+1})},
\end{align}
the following properties are satisfied:
\begin{enumerate}
\item $\Gamma(G_n(R),p^r)/\Gamma(G_n(R),p^{r+1})\cong\bigoplus_{n^2-1}{\F_p[V]}$.
\item $\Gamma(G_n(R),p^r)/\Gamma(G_n(R),p^{r+1})\cong\mathfrak{sl}_n(\F_p[V])$.
\item $\psi_r^p(A_{ij,k,r})=A_{ij,k,r+1}$, where the matrices $\{A_{ij,k,r}\}$ are as defined in Corollary \ref{SL.generators}.
\item $\gr_{*}(\Gamma(G_n(R),p))$ is generated as a restricted Lie algebra by the matrices $\{A_{ij,k,1}\}$.
\item The Lie bracket is defined on the generators $\{A_{ij,k,r}\}$ as follows:
\footnotesize
$$[A_{ij,q_1,r},A_{kl,q_2,s}]=\begin{cases} A_{kj,q_1q_2,r+s}^{-1}                  &\text{if $i=l$, $j\neq k$, $i\neq j$, $k\neq l$}\\
                                            A_{ii,q_1q_2,r+s}A_{jj,q_1q_2,r+s}^{-1} &\text{if $i=l$, $j=k$, $i\neq j$, $k\neq l$}\\
                                            A_{il,q_1q_2,r+s}                       &\text{if $i\neq l$, $j=k$, $i\neq j$, $k\neq l$}\\
                                            A_{in,q_1q_2,r+s}^{-2}                  &\text{if $i=k=l$, $i\neq j$, $j=n$}\\
                                            A_{ij,q_1q_2,r+s}^{-1}                  &\text{if $i=k=l$, $i\neq j$, $j\neq n$, $i\neq n$}\\
                                            A_{nj,q_1q_2,r+s}^2                     &\text{if $j=k=l$, $i\neq j$, $i=n$}\\
                                            A_{ij,q_1q_2,r+s}                       &\text{if $j=k=l$, $i\neq j$, $i\neq n$, $j\neq n$}\\
                                            A_{in,q_1q_2,r+s}^{-1}                  &\text{if $k=l$, $i\neq k$, $j\neq k$, $i\neq j$, $j=n$}\\
                                            A_{nj,q_1q_2,r+s}                       &\text{if $k=l$, $i\neq k$, $j\neq k$, $i\neq j$, $i=n$}\\
                                            A_{ni,q_1q_2,r+s}^{-2}                  &\text{if $i=j=l$, $k\neq l$, $k=n$}\\
                                            A_{ki,q_1q_2,r+s}^{-1}                  &\text{if $i=j=l$, $k\neq l$, $k\neq n$, $l\neq n$}\\
                                            A_{in,q_1q_2,r+s}^2                     &\text{if $i=j=k$, $k\neq l$, $l=n$}\\
                                            A_{il,q_1q_2,r+s}                       &\text{if $i=j=k$, $k\neq l$, $k\neq n$, $l\neq n$}\\
                                            A_{nl,q_1q_2,r+s}^{-1}                  &\text{if $i=j$, $k\neq l$, $i\neq l$, $i\neq k$, $k=n$}\\
                                            A_{kn,q_1q_2,r+s}                       &\text{if $i=j$, $k\neq l$, $i\neq l$, $i\neq k$, $l=n$}\\
                                            1                                       &\text{otherwise.}
\end{cases}$$
\normalsize
\end{enumerate}
\end{thm}
                        
\begin{proof} The first four parts of the theorem follow immediately from Lemma \ref{SL.quotients}, Theorem \ref{general.quotients.iso.to.sl.lie.algebra}, Proposition \ref{SL.p.power}, and Corollary \ref{SL.generators}, respectively.  By Lemma \ref{SL.commutator}, we have the following:
\small
$$[A_{ij,q_1,r},A_{kl,q_2,s}]=\begin{cases} 1+p^{r+s}v_{q_1}v_{q_2}[e_{ij},e_{kl}]                                   &\text{if $i\neq j$, $k\neq l$}\\
                                            1+p^{r+s}v_{q_1}v_{q_2}([e_{ij},e_{kk}]+[e_{nn},e_{ij}])                 &\text{if $i\neq j$, $k=l$}\\
                                            1+p^{r+s}v_{q_1}v_{q_2}([e_{ii},e_{kl}]+[e_{kl},e_{nn}])                 &\text{if $i=j$, $k\neq l$}\\
                                            1+p^{r+s}v_{q_1}v_{q_2}([e_{ii},e_{kk}]+[e_{nn},e_{ii}]+[e_{kk},e_{nn}]) &\text{if $i=j$, $k=l$.}
\end{cases}$$
\normalsize
Using Lemma \ref{eij.bracket}, the reader can verify the relations listed in the statement of the theorem.
\end{proof}

\section{Cohomology of congruence subgroups for $SL_2(\Z[t])$}

Suppose $n=2$ and $R=\Z[t]$, so that $G_n(R)=G_2(\Z[t])=SL_2(\Z[t])$, where $\Z[t]$ is a polynomial ring in one indeterminate.  For fixed $r\geq1$ and $i\geq0$, consider the matrix
\begin{align}
A_{12,i,r}=\begin{pmatrix}1&p^rt^i\\0&1\end{pmatrix}\in\Gamma(G_2(\Z[t]),p^r).
\end{align}
The collection $\{A_{12,i,r}\}_{i\geq0}$ is a family of pairwise commutative matrices, each of infinite order in $\Gamma(G_2(\Z[t]),p^r)$.  This collection generates a copy of $\bigoplus_{i\geq0}\Z\subseteq\Gamma(G_2(\Z[t]),p^r)$.  Thus, we have the diagram
\begin{figure}[h!]
\begin{diagram}
\bigoplus_{i\geq0}\Z & \rInto^{\phantom{12345}i\phantom{12345}} & \Gamma(G_2(\Z[t]),p^r) & \rTo^{\phantom{12345}\pi\phantom{12345}} & \Gamma(G_2(\Z[t]),p^r)/\Gamma(G_2(\Z[t]),p^{r+1})\\
                     & \rdTo(4,2)_{f}                       &                    &                                          & \dTo_{\Phi^{-1}_2}\\
                     &                                      &                    &                                          & \F_p[t]
\end{diagram}
\caption{}
\end{figure}\\
where $i$ is the inclusion map, $\pi$ is the natural quotient map, $\Phi^{-1}_2$ is the projection onto the second coordinate under the isomorphism $\Phi^{-1}$ defined in the proof of Lemma \ref{SL.quotients}, and $f$ is the composite.  Notice that $f(A_{12,i,r})=t^i$ under this mapping.  We now consider the cases $p=2$ and $p>2$ separately.

\vskip .1in

\noindent\underline{$p=2$:} Recall that
\begin{align}
H^*(\F_2[t];\F_2)\cong H^*(\oplus_{i\geq0}\F_2;\F_2)\cong\F_2[x_0,x_1,\ldots],
\end{align}
a polynomial ring in an infinite number of indeterminates of degree one, where each $x_i$ is dual to the image of $A_{12,i,r}$ in $\Gamma(G_2(\Z[t]),2^r)/\Gamma(G_2(\Z[t]),2^{r+1})$.  Also,
\begin{align}
H^*(\oplus_{i\geq0}\Z;\F_2)\cong\Lambda^*(y_0,y_1,\ldots;\F_2),
\end{align}
an exterior algebra in an infinite number of indeterminates of degree one over $\F_2$, where each $y_i$ is dual to $A_{12,i,r}$ in $\Gamma(G_2(\Z[t]),2^r)$.

As illustrated in Figure \ref{cohomology.sl2z(t).p=2}, there are induced maps in cohomology 
\begin{figure}[h!]
\begin{diagram}
H^j(\F_2[t];\F_2) & \rTo^{\phantom{12}(\Phi^{-1}_2\circ\pi)^*\phantom{12}} & H^j(\Gamma(G_2(\Z[t]),2^r);\F_2) & \rTo^{\phantom{123}i^*\phantom{123}} & H^j(\oplus_{i\geq0}\Z;\F_2)\\
                  & \rdTo(4,2)_{f^*}     &                              &            & \dEq_{\phantom{12345}}\\
                  &                           &                              &            & H^j(\oplus_{i\geq0}\Z;\F_2)
\end{diagram}
\caption{}\label{cohomology.sl2z(t).p=2}
\end{figure}
for which $f^*(x_i)=y_i$ for $j=1$.  This implies that $f^*$ is an epimorphism for $j\geq1$, forcing $i^*$ to be an epimorphism for $j\geq1$ as well.  In this case, $H^j(\bigoplus_{i\geq0}\Z;\F_2)$ is a quotient of $H^j(\Gamma(G_2(\Z[t]),2^r);\F_2)$.  Since  $H^j(\bigoplus_{i\geq0}\Z;\F_2)$ is not finitely generated for $j\geq1$, neither is $H^j(\Gamma(G_2(\Z[t]),2^r);\F_2)$.

\vskip .1in

\noindent\underline{$p>2$:}  Recall that
\begin{align}
H^*(\F_p[t];\F_p)\cong H^*(\oplus_{i\geq0}\F_p;\F_p)\cong\Lambda^*(x_0,x_1,\ldots;\F_p)\otimes\F_p[\beta x_0,\beta x_1,\ldots],
\end{align}
the tensor product of an exterior algebra with a polynomial ring.  Here, each $x_i$ has degree one, and $\begin{CD}\beta:H^1(\F_p;\F_p)@>>> H^2(\F_p;\F_p)\end{CD}$ is the Bockstein homomorphism.  As above, each $x_i$ is dual to the image of $A_{12,i,r}$ in $\Gamma(G_2(\Z[t]),p^r)/\Gamma(G_2(\Z[t]),p^{r+1})$.  Also,
\begin{align}
H^*(\oplus_{i\geq0}\Z;\F_p)\cong\Lambda^*(y_0,y_1,\ldots;\F_p),
\end{align}
an exterior algebra in an infinite number of indeterminates of degree one over $\F_p$, where each $y_i$ is dual to $A_{12,i,r}$ in $\Gamma(G_2(\Z[t]),p^r)$.

There are induced maps in cohomology
\begin{figure}[h!]
\begin{diagram}
H^j(\F_p[t];\F_p) & \rTo^{\phantom{12}(\Phi^{-1}_2\circ\pi)^*\phantom{12}} & H^j(\Gamma(G_2(\Z[t]),p^r);\F_p) & \rTo^{\phantom{123}i^*\phantom{123}} & H^j(\oplus_{i\geq0}\Z;\F_p)\\
                  & \rdTo(4,2)_{f^*}     &                              &            & \dEq_{\phantom{12345}}\\
                  &                           &                              &            & H^j(\oplus_{i\geq0}\Z;\F_p)
\end{diagram}
\caption{}
\end{figure}\\
for which $f^*(x_i)=y_i$ for $j=1$.  This implies that $f^*$ is an epimorphism for $j\geq1$, forcing $i^*$ to be an epimorphism for $j\geq1$ as well.  In this case, $H^j(\bigoplus_{i\geq0}\Z;\F_p)$ is a quotient of $H^j(\Gamma(G_2(\Z[t]),p^r);\F_p)$.  Since  $H^j(\bigoplus_{i\geq0}\Z;\F_p)$ is not finitely generated for $j\geq1$, neither is $H^j(\Gamma(G_2(\Z[t]),p^r);\F_p)$.  Thus, we have proved the following theorem.

\begin{thm}\label{infinitely.generated.cohomology}
Suppose $n=2$ and $R=\Z[t]$.  For all $j,r\geq1$ and for all primes $p$,
\begin{align}
H^j(\Gamma(G_2(\Z[t]),p^r);\F_p)
\end{align}
is infinitely generated.
\end{thm}

\section{Cohomology of congruence subgroups for $SL_2(\Z[i])$}

Suppose $n=2$ and $R=\Z[i]$, so that $G_n(R)=G_2(\Z[i])=SL_2(\Z[i])$, where $\Z[i]$ is the Gaussian integers.  Let $e_{ij}$ denote the $2\times2$ matrix with a 1 in the $(i,j)$ position and zeros elsewhere.  For fixed $r\geq1$, $\epsilon\in\{0,1\}$, and $i\neq j$, consider the matrices $A_{ij,\epsilon,r}=1+p^ri^{\epsilon}e_{ij}\in\Gamma(G_2(\Z[i]),p^r)$.  Each of these four matrices has infinite order in $\Gamma(G_2(\Z[i]),p^r)$.  Notice that $A_{12,0,r}$ and $A_{12,1,r}$ commute, so they generate a copy of $\Z\oplus\Z\subseteq\Gamma(G_2(\Z[i]),p^r)$.  The same is true for the matrices $A_{21,0,r}$ and $A_{21,1,r}$.  Thus, we have embeddings $\begin{CD}i_1:\Z\oplus\Z@>>>\Gamma(G_2(\Z[i]),p^r)\end{CD}$ and $\begin{CD}i_2:\Z\oplus\Z@>>>\Gamma(G_2(\Z[i]),p^r)\end{CD}$.  These embeddings give the following diagram:
\small
\begin{figure}[H]
\begin{diagram}
(\Z\oplus\Z)*(\Z\oplus\Z) & \rTo^{\phantom{1}i_1*i_2\phantom{1}} & \Gamma(G_2(\Z[i]),p^r) & \rTo^{\phantom{1234}\pi\phantom{1234}} & \Gamma(G_2(\Z[i]),p^r)/\Gamma(G_2(\Z[i]),p^{r+1})\\
                     & \rdTo(4,2)_{f}                       &                    &                                          & \dTo_{\Phi^{-1}}\\
                     &                                      &                    &                                          & \oplus_6{\F_p}
\end{diagram}
\caption{}
\end{figure}
\normalsize
\noindent Here, $\pi$ is the natural quotient map, $\Phi$ is the isomorphism defined in the proof of Lemma \ref{SL.quotients}, and $f$ is the composite.  We now consider the cases $p=2$ and $p>2$ separately.

\vskip .1in
	
\noindent\underline{$p=2$:}  Recall that
\begin{align}
H^*(\oplus_6\F_2;\F_2)\cong\F_2[x_1,x_2,y_1,y_2,z_1,z_2],
\end{align}
a polynomial ring in six indeterminates of degree one, where $x_1$ is dual to the image of $A_{12,0,r}$, $x_2$ is dual to the image of $A_{12,1,r}$, $y_1$ is dual to the image of $A_{21,0,r}$, and $y_2$ is dual to the image of $A_{21,1,r}$ in $\Gamma(G_2(\Z[i]),2^r)/\Gamma(G_2(\Z[i]),2^{r+1})$.  Also,
\begin{align}
H^*((\Z\oplus\Z)*(\Z\oplus\Z);\F_2)\cong\Lambda^*(u_1,u_2;\F_2)\oplus\Lambda^*(w_1,w_2;\F_2),
\end{align}
where $u_1$ is dual to $A_{12,0,r}$, $u_2$ is dual to $A_{12,1,r}$, $w_1$ is dual to $A_{21,0,r}$, and $w_2$ is dual to $A_{21,1,r}$ in $\Gamma(G_2(\Z[i]),2^r)$.  All of $u_1$, $u_2$, $w_1$, and $w_2$ have degree one.  In particular, notice that
\begin{align}
H^2((\Z\oplus\Z)*(\Z\oplus\Z);\F_2)\cong\F_2\oplus\F_2,
\end{align}
generated by the classes $u_1u_2$ and $w_1w_2$.

There are induced maps in cohomology
\small
\begin{figure}[h!]
\begin{diagram}
H^j(\oplus_6\F_2;\F_2) & \rTo^{\phantom{1}(\Phi^{-1}\circ\pi)^*\phantom{1}} & H^j(\Gamma(G_2(\Z[i]),2^r);\F_2) & \rTo^{(i_1*i_2)^*} & H^j((\Z\oplus\Z)*(\Z\oplus\Z);\F_2)\\
                  & \rdTo(4,2)_{f^*}     &                              &            & \dEq_{\phantom{12345}}\\
                  &                      &                              &            & H^j((\Z\oplus\Z)*(\Z\oplus\Z);\F_2)
\end{diagram}
\caption{}
\end{figure}\\
\normalsize
for which $f^*(x_i)=u_i$ and $f^*(y_i)=w_i$ for $j=1$.  Thus, we see that $f^*$ is an epimorphism for $j=2$, forcing $(i_1*i_2)^*$ to be an epimorphism for $j=2$ as well.  As a result, $H^2((\Z\oplus\Z)*(\Z\oplus\Z);\F_2)$ is a quotient of $H^2(\Gamma(G_2(\Z[i]),2^r);\F_2)$, so that $\F_2\oplus\F_2\subseteq H^2(\Gamma(G_2(\Z[i]),2^r);\F_2)$.

\vskip .1in

\noindent\underline{$p>2$:}  Recall that
\begin{align}
H^*(\oplus_6\F_p;\F_p)\cong\Lambda^*(x_1,x_2,y_1,y_2,z_1,z_2;\F_p)\otimes\F_p[\beta x_1,\beta x_2,\beta y_1,\beta y_2,\beta z_1,\beta z_2],
\end{align}
the tensor product of an exterior algebra with a polynomial ring.  Here, each $x_i$, $y_i$, and $z_i$ is of degree one and $\begin{CD}\beta:H^1(\F_p;\F_p)@>>>H^2(\F_p;\F_p)\end{CD}$ is the Bockstein homomorphism.  Similar to above, $x_1$ is dual to the image of $A_{12,0,r}$, $x_2$ is dual to the image of $A_{12,1,r}$, $y_1$ is dual to the image of $A_{21,0,r}$, and $y_2$ is dual to the image of $A_{21,1,r}$ in $\Gamma(G_2(\Z[i]),p^r)/\Gamma(G_2(\Z[i]),p^{r+1})$.  Also,
\begin{align}
H^*((\Z\oplus\Z)*(\Z\oplus\Z);\F_p)\cong\Lambda^*(u_1,u_2;\F_p)\oplus\Lambda^*(w_1,w_2;\F_p),
\end{align}
where $u_1$ is dual to $A_{12,0,r}$, $u_2$ is dual to $A_{12,1,r}$, $w_1$ is dual to $A_{21,0,r}$, and $w_2$ is dual to $A_{21,1,r}$ in $\Gamma(G_2(\Z[i]),p^r)$.  All of $u_1$, $u_2$, $w_1$, and $w_2$ have degree one.  In particular, notice that
\begin{align}
H^2((\Z\oplus\Z)*(\Z\oplus\Z);\F_p)\cong\F_p\oplus\F_p,
\end{align}
generated by the classes $u_1u_2$ and $w_1w_2$.

As illustrated in Figure \ref{cohom.diag.gaussian}, there are induced maps in cohomology 
\small
\begin{figure}[h!]
\begin{diagram}
H^j(\oplus_6\F_p;\F_p) & \rTo^{\phantom{1}(\Phi^{-1}\circ\pi)^*\phantom{1}} & H^j(\Gamma(G_2(\Z[i]),p^r);\F_p) & \rTo^{(i_1*i_2)^*} & H^j((\Z\oplus\Z)*(\Z\oplus\Z);\F_p)\\
                  & \rdTo(4,2)_{f^*}     &                              &            & \dEq_{\phantom{12345}}\\
                  &                      &                              &            & H^j((\Z\oplus\Z)*(\Z\oplus\Z);\F_p)
\end{diagram}
\caption{}\label{cohom.diag.gaussian}
\end{figure}
\normalsize
for which $f^*(x_i)=u_i$ and $f^*(y_i)=w_i$ for $j=1$.  Thus, we see that $f^*$ is an epimorphism for $j=2$, forcing $(i_1*i_2)^*$ to be an epimorphism for $j=2$ as well.  As a result, $H^2((\Z\oplus\Z)*(\Z\oplus\Z);\F_p)$ is a quotient of $H^2(\Gamma(G_2(\Z[i]),p^r);\F_p)$, so that $\F_p\oplus\F_p\subseteq H^2(\Gamma(G_2(\Z[i]),p^r);\F_p)$.  We have now proved the following theorem.

\begin{thm}\label{bound.gaussian.cohomology}
Suppose $n=2$ and $R=\Z[i]$.  For all $r\geq1$ and for all primes $p$,
\begin{align}
H^2(\Gamma(G_2(\Z[i]),p^r);\F_p)\supseteq\F_p\oplus\F_p.
\end{align}
\end{thm}

\section{Computing $d^2$ for central extensions of congruence subgroups}

In this section, we use the convention $\Gamma_r=\Gamma(G_n(R),p^r)$.  Let $F[x,y]$ be the free group on two letters, and let $K$ denote the kernel of the abelianization map $\begin{CD}F[x,y]@>>>\Z\oplus\Z\end{CD}$.

Suppose $r,s,l\geq1$.  Specifying elements $\theta(x),\theta(y)\in\Gamma_r/\Gamma_{r+s}$ defines a unique homomorphism
\begin{align}
\begin{CD}\theta:F[x,y]@>>>\Gamma_r/\Gamma_{r+s}.
\end{CD}
\end{align}
We are interested in situations in which two conditions are satisfied:
\begin{align}\label{congruence.extension}
\begin{CD}
1@>>>\Gamma_{r+s-l}/\Gamma_{r+s}@>>>\Gamma_r/\Gamma_{r+s}@>>>\Gamma_r/\Gamma_{r+s-l}@>>>1
\end{CD}
\end{align}
is central; and defining a map $\theta$ as above induces a morphism of extensions
\begin{figure}[h!]
$\begin{CD}
  1@>>>K@>>>F[x,y]@>>>\Z\oplus\Z @>>>1 \\
  @. @VV{\theta}V @VV{\theta}V @VV{\theta}V @. \\
  1@>>>\Gamma_{r+s-l}/\Gamma_{r+s}@>>>\Gamma_{r}/\Gamma_{r+s}@>>>\Gamma_{r}/\Gamma_{r+s-l}@>>>1.
\end{CD}$
\caption{}
\end{figure}\\
It is straightforward to determine when the former condition is satisfied.

\begin{lem}\label{SL.central.criteria}  
For $r,s,l\geq1$, extension (\ref{congruence.extension}) is central if and only if $r\geq l$.
\end{lem}

\begin{proof}
In the proof of Theorem \ref{SL.s.stage.quotients} (which is included in Section \ref{SL.s.stage.quotients.proof}), we show that $\Gamma(G_n(R),p^r)/\Gamma(G_n(R),p^{r+s})$ is generated by the images of the matrices $\{A_{ij,k,r}\}$ and that $\Gamma(G_n(R),p^{r+s-l})/\Gamma(G_n(R),p^{r+s})$ is generated by the images of the matrices $\{A_{ij,k,r+s-l}\}$.  A straightforward calculation will verify that the following relation is satisfied, where $[g,h]=g^{-1}h^{-1}gh$ denotes the group commutator:
\begin{align}
[A_{i_1j_1,k_1,r+s-l},A_{i_2j_2,k_2,r}]=1+p^{2r+s-l}v_{k_1}v_{k_2}(e_{i_1j_1}e_{i_2j_2}-e_{i_2j_2}e_{i_1j_1}).
\end{align}
From this, we see that the commutator is trivial in $\Gamma_r/\Gamma_{r+s}$ if and only if $2r+s-l\geq{r+s}$.
\end{proof}

If $l=1$ and $r\geq s-1$, one can check that extension (\ref{congruence.extension}) also satisfies the second condition mentioned above.  Namely, defining a map $\begin{CD}\theta:F[x,y]@>>>\Gamma_r/\Gamma_{r+s}\end{CD}$ induces a morphism of extensions
\begin{figure}[h!]
$\begin{CD}
  1@>>>K@>>>F[x,y]@>>>\Z\oplus\Z @>>>1 \\
  @. @VV{\theta}V @VV{\theta}V @VV{\theta}V @. \\
  1@>>>\Gamma_{r+s-1}/\Gamma_{r+s}@>>>\Gamma_{r}/\Gamma_{r+s}@>>>\Gamma_{r}/\Gamma_{r+s-1}@>>>1. \\
\end{CD}$
\caption{}\label{induced.map.of.extensions}
\end{figure}

Let $\chi$ denote the generator of $H_2(\Z\oplus\Z;\Z)\cong\Z$.  Working out the local coefficient system in the Lyndon-Hochschild-Serre spectral sequence for the extension
\begin{align}
\begin{CD}
1@>>>K@>>>F[x,y]@>>>\Z\oplus\Z@>>>1,
\end{CD}
\end{align}
one can check that $H_0(\Z\oplus\Z;H_1(K))\cong\Z$.  Furthermore, this group is generated by the class of the element $[x,y]\in K$.  Since it is known that $H_2(F[x,y];\Z)=0$, it must be the case that the differential $d^2_{2,0}$ is an isomorphism.  In other words, $d^2_{2,0}(\chi)=[x,y]$.

The induced map $\theta_*$ on homology induces a map of spectral sequences associated to the two extensions in Figure \ref{induced.map.of.extensions}.  Using naturality and the fact that $d^2_{2,0}(\chi)=[x,y]$, we can write
\begin{align}
d^2_{2,0}(\theta_*(\chi))=\theta_*(d^2_{2,0}(\chi))=\theta_*([x,y])=[\theta_*(x),\theta_*(y)].
\end{align}
Thus, we have demonstrated the following formula.

\begin{thm}\label{congruence.differential}  
Suppose $r,s\geq1$, $r\geq s-1$, and a map $\begin{CD}\theta:F[x,y]@>>>\Gamma_r/\Gamma_{r+s}\end{CD}$ is given.  Then $\theta$ induces a morphism of extensions as in Figure \ref{induced.map.of.extensions}, and
\begin{align}
d^2_{2,0}(\theta_*(\chi))=[\theta_*(x),\theta_*(y)]
\end{align}
in the LHS spectral sequence for the extension
\begin{align}\label{congruence.extension.1}
\begin{CD}
1@>>>\Gamma_{r+s-1}/\Gamma_{r+s}@>>>\Gamma_{r}/\Gamma_{r+s}@>>>\Gamma_{r}/\Gamma_{r+s-1}@>>>1.
\end{CD}
\end{align}
Here, $\chi$ denotes the generator of $H_2(\Z\oplus\Z;\Z)\cong\Z$.
\end{thm}

Notice that $[\theta_*(x),\theta_*(y)]$ can be calculated using the Lie algebra structure for $\gr_{*}(\Gamma(G_n(R),p))$ given in Theorem \ref{SL(n,Z[V]).lie.algebra}.  This will allow us to compute $d^2_{2,0}$ for many examples involving central extensions of congruence subgroups.  This will be discussed in more detail in a sequel.

\vskip .1in

\noindent\textit{Remark.}  The reader should notice that the formula for the differential $d^2_{2,0}$ in Theorem \ref{congruence.differential} provides information that is dual, in a sense, to the extension class characterizing extension (\ref{congruence.extension.1}).

\section{Proof of Theorem \ref{general.quotients.iso.to.sl.lie.algebra}}\label{general.quotients.iso.to.sl.lie.algebra.proof}

We adapt the method of the proof of Theorem \ref{homology.level.p} as found in \cite{knudson.book}.  Suppose that $X\in\Gamma(G_n(R),p^r)$.  We can write $X=1+p^rA$ where 1 denotes the $n\times n$ identity matrix and $A\in\Mat_n(R)$.  Define the map 
\begin{align}
\begin{CD}
\varphi_r:\Gamma(G_n(R),p^r)@>>>\mathfrak{sl}_n(\F_p[V])
\end{CD}
\end{align}
by
\begin{align}
\varphi_r(X)=A\mod p.
\end{align}
To show that the image of $\varphi_r$ lies in $\mathfrak{sl}_n(\F_p[V])$, we make use of the fact that $\Gamma(G_n(R),p^r)\subseteq SL_n(R)$.  Thus,
	\begin{align*}
1       & =     \det(X)\\
        & =     \det(1+p^rA)\\
        &\equiv (1+p^rtr(A))\mod p^{r+1},
\end{align*}
where $tr(A)$ denotes the trace of $A$.  Thus, it must be the case that $tr(A)\equiv0\mod p$.

One can check that $\varphi_r$ is a homomorphism directly:
\begin{align*}
\varphi_r((1+p^rA)(1+p^rB))&=      \varphi_r(1+p^r(A+B+p^rAB))\\
                           &=      (A+B+p^rAB)\mod p\\
                           &\equiv (A+B)\mod p\\
                           &=      \varphi_r(1+p^rA)+\varphi_r(1+p^rB).
\end{align*}

To show that $\varphi_r$ is surjective, we define a basis for $\mathfrak{sl}_n(\F_p[V])$ and show that $\varphi_r$ surjects onto this basis.  In the proof of Corollary \ref{homology.level.p.general}, we defined $e_{ij}$ to be the $n\times n$ matrix with a 1 in the $(i,j)$ position and zeros elsewhere.  A quick check will verify that a basis for $\mathfrak{sl}_n(\F_p[V])$ is given by $\{v_ke_{ij}\}_{i\neq j}\cup\{v_k(e_{ii}-e_{nn})\}_{i=1}^{n-1}$ for $k\in I$.  Further,
\begin{align}
\varphi_r(1+p^rv_ke_{ij})=v_ke_{ij}
\end{align}
and
\begin{align}
\varphi_r(1+p^rv_k(e_{ii}+e_{in}-e_{ni}-e_{nn}))=v_k(e_{ii}+e_{in}-e_{ni}-e_{nn}),
\end{align}
so that $\varphi_r$ hits all of the basis elements in $\mathfrak{sl}_n(\F_p[V])$.

Finally, notice that if $X\in\Gamma(G_n(R),p^{r+1})$, we can write $X=1+p^{r+1}A=1+p^r(pA)$ for some $A\in\Mat_n(R)$.  So $\varphi_r(X)=0$.  Conversely, if $\varphi_r(1+p^rA)=0$, it must be the case that $A=pB$ for some $B\in\Mat_n(R)$.  Thus, $1+p^rA=1+p^{r+1}B\in\Gamma(G_n(R),p^{r+1})$.   So $\ker(\varphi_r)=\Gamma(G_n(R),p^{r+1})$.  This completes the proof of the theorem.

\section{Proof of Theorem \ref{SL.s.stage.quotients}}\label{SL.s.stage.quotients.proof}

The proof is by induction on $s$.  The case $s=1$ is Lemma \ref{SL.quotients}.  Suppose that $s-1$ and smaller cases have been proved and that $|V|<\infty$.  Consider the group extension
\begin{align}
\begin{CD}
1@>>>\Gamma_{r+s-1}/\Gamma_{r+s}@>>>\Gamma_r/\Gamma_{r+s}@>>>\Gamma_r/\Gamma_{r+s-1}@>>>1,
\end{CD}
\end{align}
where we again use the convention $\Gamma_l=\Gamma(G_n(R),p^l)$.  By Lemma \ref{SL.quotients},
\begin{align}
\Gamma_{r+s-1}/\Gamma_{r+s}\cong\bigoplus_{n^2-1}{\Z/p\Z[V]},
\end{align}
so that $|\Gamma_{r+s-1}/\Gamma_{r+s}|=p^{(n^2-1)\cdot|V|}$.  By induction, 
\begin{align}
\Gamma_r/\Gamma_{r+s-1}\cong\bigoplus_{n^2-1}{\Z/p^{s-1}\Z[V]},
\end{align}
so that $|\Gamma_r/\Gamma_{r+s-1}|=p^{(s-1)\cdot(n^2-1)\cdot|V|}$.  It follows that $|\Gamma_r/\Gamma_{r+s}|=p^{(n^2-1)\cdot|V|}\cdot p^{(s-1)\cdot(n^2-1)\cdot|V|}=p^{s\cdot(n^2-1)\cdot|V|}$.

Any element of $\bigoplus_{n^2-1}{\Z/{p^s}\Z[V]}$ can be written as
\begin{align}
(a_{11},a_{12},\ldots,a_{1n},a_{21},\ldots,a_{2n},\ldots,a_{n1},\ldots,a_{n,n-1}),
\end{align}
where $a_{ij}\in\Z/{p^s}\Z[V]$.  Let $\delta_{ij}\in\bigoplus_{n^2-1}{\Z/{p^s}\Z[V]}$ be the element with a 1 in the the $(i,j)$ position and zeros elsewhere.  Then $\{v_k\delta_{ij}\}$ for $1\leq i,j\leq n$, $i+j<2n$, and $V=\{v_k\}_{k\in I}$ is a basis for $\bigoplus_{n^2-1}{\Z/{p^s}\Z[V]}$.  Consider the map 
\begin{align}
\begin{CD}
\Phi_{r,s}:\bigoplus_{n^2-1}{\Z/{p^s}\Z[V]}@>>>{\Gamma_r/\Gamma_{r+s}}
\end{CD}
\end{align}
that sends $v_k\delta_{ij}\mapsto\overline{A}_{ij,k,r}$, where $\overline{A}_{ij,k,r}$ represents the equivalence class of $A_{ij,k,r}$ in $\Gamma_r/\Gamma_{r+s}$, and the family $\{A_{ij,k,r}\}$ are the generators of $\Gamma_r/\Gamma_{r+1}$ as defined in Corollary \ref{SL.generators}.

\vskip .1in\noindent\textit{Claim.}  The group $\Gamma_r/\Gamma_{r+s}$ is generated by the equivalence classes $\{\overline{A}_{ij,k,r}\}$ for $1\leq i,j\leq n$, $i+j<2n$, and $k\in I$.

\vskip .1in\noindent\textit{Proof of Claim.}  We again proceed by induction on $s$.  The case $s=1$ is clear.  Suppose that $s-1$ and smaller cases are satisfied.  By induction, $\Gamma_r/\Gamma_{r+s-1}$ is generated by the equivalence classes of the matrices $\{A_{ij,k,r}\}$.  According to Corollary \ref{SL.generators}, $\Gamma_{r+s-1}/\Gamma_{r+s}$ is generated by the matrices $\{A_{ij,k,r+s-1}\}$.  But for each $i$, $j$, and $k$, 
\begin{align}
A_{ij,k,r+s-1}=(\psi_{r+s-2}^p\circ\psi_{r+s-3}^p\circ\cdots\circ\psi_{r+1}^p\circ\psi_r^p)(A_{ij,k,r})=A_{ij,k,r}^{p^{s-1}}.
\end{align}
Thus, $\Gamma_r/\Gamma_{r+s}$ is generated by the equivalence classes of the matrices $\{A_{ij,k,r}\}$, which completes the proof of the claim.

Once we have this, the reader can check that since $r\geq s\geq2$, the generators of $\Gamma_r/\Gamma_{r+s}$ commute, so that $\Gamma_r/\Gamma_{r+s}$ is abelian.  So $\Phi_{r,s}$ is a homomorphism of abelian groups, which is clearly surjective since $\Gamma_r/\Gamma_{r+s}$ is generated by the equivalence classes of the matrices $\{A_{ij,k,r}\}$.

Finally, since $|\bigoplus_{n^2-1}{\Z/{p^s}\Z[V]}|=p^{s\cdot(n^2-1)\cdot|V|}=|\Gamma_r/\Gamma_{r+s}|$, it must also be the case that $\Phi_{r,s}$ is injective, making it an isomorphism.  This completes the proof of Theorem \ref{SL.s.stage.quotients}.

\section{Proof of Theorem \ref{main.theorem}}\label{proof.of.main.theorem}

Recall the map 
\begin{align}
\begin{CD}
\varphi_r:\Gamma(G_n(R),p^r)@>>>\mathfrak{sl}_n(\F_p[V])
\end{CD}
\end{align}
defined in Section \ref{general.quotients.iso.to.sl.lie.algebra.proof}.  This induces a well-defined map on the filtration quotients
\begin{align}
\begin{CD}
\overline{\varphi}_r:\Gamma(G_n(R),p^r)/\Gamma(G_n(R),p^{r+1})@>>>\mathfrak{sl}_n(\F_p[V]).
\end{CD}
\end{align}
In fact, since $\varphi_r$ is an epimorphism with kernel $\Gamma(G_n(R),p^{r+1})$ (as shown in Section \ref{general.quotients.iso.to.sl.lie.algebra.proof}), $\overline{\varphi}_r$ is an isomorphism.

Recall that $I$ denotes the kernel of the evaluation map
$\begin{CD}
  \F_p[t]@>{t=0}>>\F_p.
\end{CD}$
Define a map
\begin{align}
\begin{CD}
\overline{\varphi}_r\otimes t^r:\Gamma(G_n(R),p^r)/\Gamma(G_n(R),p^{r+1})@>>>\mathfrak{sl}_n(\F_p[V])\otimes_{\F_p}I
\end{CD}
\end{align}
by
\begin{align}
(\overline{\varphi}_r\otimes t^r)(X_r)=\overline{\varphi}_r(X_r)\otimes t^r
\end{align}
for $X_r\in\Gamma(G_n(R),p^r)/\Gamma(G_n(R),p^{r+1})$.  We can assemble these maps for $r\geq1$ to obtain a map 
\begin{align}
\begin{CD}
\varphi:\gr_{*}(\Gamma(G_n(R),p))@>>>\mathfrak{sl}_n(\F_p[V])\otimes_{\F_p}I
\end{CD}
\end{align}
given by
\begin{align}
\varphi(X)=\sum_{r\geq1}(\overline{\varphi}_r\otimes t^r)(X_r).
\end{align}
Here, we have made use of the fact that
\begin{align}
\gr_{*}(\Gamma(G_n(R),p))=\bigoplus_{r\geq1}\Gamma(G_n(R),p^r)/\Gamma(G_n(R),p^{r+1})
\end{align}
and written $X=(X_1,X_2,\ldots)$, where $X_r\in\Gamma(G_n(R),p^r)/\Gamma(G_n(R),p^{r+1})$.  We claim that $\varphi$ is the required Lie algebra isomorphism.

We first verify that $\varphi$ is a homomorphism.  Suppose $X,Y\in\gr_{*}(\Gamma(G_n(R),p))$.  Write $X=(X_1,X_2,\ldots)$ and $Y=(Y_1,Y_2,\ldots)$, where
\begin{align}
X_r,Y_r\in\Gamma(G_n(R),p^r)/\Gamma(G_n(R),p^{r+1})
\end{align}
for $r\geq1$.  Since each $\overline{\varphi}_r$ is a homomorphism, we have the following:
\begin{align*}
\varphi(X+Y)& = \sum_{r\geq1}{(\overline{\varphi}_r\otimes t^r)(X_rY_r)}\\[.1em]
            & = \sum_{r\geq1}{(\overline{\varphi}_r(X_rY_r)\otimes t^r)}\\[.1em]
            & = \sum_{r\geq1}{([\overline{\varphi}_r(X_r)+\overline{\varphi}_r(Y_r)]\otimes t^r)}\\[.1em]
            & = \sum_{r\geq1}{([\overline{\varphi}_r(X_r)\otimes t^r]+[\overline{\varphi}_r(Y_r)\otimes t^r])}\\[.1em]
            & = \sum_{r\geq1}{(\overline{\varphi}_r(X_r)\otimes t^r)}+\sum_{r\geq1}{(\overline{\varphi}_r(Y_r)\otimes t^r)}\\[.1em]
            & = \sum_{r\geq1}{(\overline{\varphi}_r\otimes t^r)(X_r)}+\sum_{r\geq1}{(\overline{\varphi}_r\otimes t^r)(Y_r)}\\[.1em]
            & = \varphi(X)+\varphi(Y).
\end{align*}
Thus, $\varphi$ is a homomorphism.

To check that $\varphi$ is surjective, we define a basis for $\mathfrak{sl}_n(\F_p[V])\otimes_{\F_p}I$ and show that $\varphi$ surjects onto this basis.  In Section \ref{general.quotients.iso.to.sl.lie.algebra.proof}, we defined a basis for $\mathfrak{sl}_n(\F_p[V])$ to be $\{v_ke_{ij}\}_{i\neq j}\cup\{v_k(e_{ii}-e_{nn})\}_{i=1}^{n-1}$ for $k\in I$.  The reader can check that $\{v_ke_{ij}\otimes t^r\}_{i\neq j}\cup\{v_k(e_{ii}-e_{nn})\otimes t^r\}_{i=1}^{n-1}$ for $k\in I$ and $r\geq1$ is a basis for $\mathfrak{sl}_n(\F_p[V])\otimes_{\F_p}I$.  It is clear that
\begin{align}
\varphi(1+p^rv_ke_{ij})=v_ke_{ij}\otimes t^r.
\end{align}
Furthermore,
\begin{align}
\varphi(1+p^rv_k(e_{ii}+e_{in}-e_{ni}-e_{nn}))=v_k(e_{ii}+e_{in}-e_{ni}-e_{nn})\otimes t^r,
\end{align}
so that $\varphi$ hits all the basis elements in $\mathfrak{sl}_n(\F_p[V])\otimes_{\F_p}I$.

Since $\overline{\varphi}_r$ is injective for all $r\geq1$ (an isomorphism, in fact), it follows that $\varphi$ is injective.

Finally, we need to show that $\varphi$ preserves the Lie bracket.  Since $[X,Y]\in\gr_{*}(\Gamma(G_n(R),p))$, we can write 
\begin{align}
[X,Y]=([X,Y]_1,[X,Y]_2,\ldots),
\end{align}
where $[X,Y]_r\in\Gamma(G_n(R),p^r)/\Gamma(G_n(R),p^{r+1})$.  Also, since
\begin{align}
X_r,Y_r\in\Gamma(G_n(R),p^r)/\Gamma(G_n(R),p^{r+1}),
\end{align}
there exist $\widehat{X}_r,\widehat{Y}_r\in\Mat_n(R)$ such that $X_r=1+p^r\widehat{X}_r$ and $Y_r=1+p^r\widehat{Y}_r$.  We have the following:
\begin{align*}
\varphi([X,Y]) & = \sum_{r\geq1}{(\overline{\varphi}_r\otimes t^r)([X,Y]_r)}\\[.1em]
               & = \sum_{r\geq1}{(\overline{\varphi}_r([X,Y]_r)\otimes t^r)}\\[.1em]
               & = \sum_{r\geq1}{(\overline{\varphi}_r(\sum_{i+j=r}{[X_i,Y_j]})\otimes t^r)}\\[.1em]
               & = \sum_{r\geq1}{(\overline{\varphi}_r(\sum_{i+j=r}{[1+p^i\widehat{X}_i,1+p^j\widehat{Y}_j]})\otimes t^r)}\\[.1em]
               & = \sum_{r\geq1}{(\sum_{i+j=r}{\overline{\varphi}_r([1+p^i\widehat{X}_i,1+p^j\widehat{Y}_j])}\otimes t^r)}\\[.1em]
               & = \sum_{r\geq1}{(\sum_{i+j=r}{((\widehat{X}_i\widehat{Y}_j-\widehat{Y}_j\widehat{X}_i)\mod p)}\otimes t^r)}\\[.1em]
               & = [\varphi(X),\varphi(Y)].
\end{align*}
Thus,
\begin{align}
\gr_{*}(\Gamma(G_n(R),p))\cong\mathfrak{sl}_n(\F_p[V])\otimes_{\F_p}I
\end{align}
as Lie algebras.  This completes the proof of the theorem.

\vskip .1in

\noindent\textit{Other Acknowledgments.}  The author would like to thank Paul Taylor for making his \textit{diagram} package available, which was used to produce several of the diagrams in this paper.


\begin{thebibliography}{20}

\bibitem{BMS}
H.~Bass, J.~Milnor, J.-P.~Serre, \emph{Solution of the congruence subgroup problem for $SL_n$ ($n\geq3$) and $Sp_{2n}$ ($n\geq2$)}, Publ. IHES \textbf{33} (1969), 421--499.

\bibitem{bigelow1}
S.~Bigelow, \emph{Braid groups are linear}, J. Amer. Math. Soc. \textbf{14} (2001), 471--486.

\bibitem{bigelow2}
S.~Bigelow, \emph{The Burau representation is not faithful for $n=5$}, Geom. Topol. \textbf{3} (1999), 397--404.

\bibitem{cohen}
F.~Cohen, M.~Condor, J.~Lopez, and S.~Prassidis, \emph{Remarks concerning Lubotzky's filtration}, Pure and Applied Mathematics Quarterly, \textbf{8} (1), pp. 79-106.  [Available at arXiv:0710.3515]

\bibitem{cl1}
D.~Cooper and D.D.~Long, \emph{On the Burau representation modulo a small prime}, Geom. Topol. Monogr. \textbf{1} (1998), 127--138.

\bibitem{cl2}
D.~Cooper and D.D.~Long, \emph{A presentation for the image of $Burau(4)\otimes\Z_2$}, Inventiones Math. \textbf{127} (1997), 535--570.

\bibitem{dusautoy}
M.~du Sautoy, J.D.~Dixon, et al., \emph{Analytic pro-p groups}, Second Edition, New York: Cambridge University Press, 1999.

\bibitem{knudson1}
K.~Knudson, \emph{Congruence subgroups and twisted cohomology of $SL(n,\F[t])$}, J. Algebra \textbf{207} (1998), 695--721.

\bibitem{knudson2}
K.~Knudson, \emph{The homology of $SL(2,\F[t,t^{-1}])$}, J. Algebra \textbf{180} (1996), 87--101.

\bibitem{knudson3}
K.~Knudson, \emph{Homology and finiteness properties for $SL(2,\Z[t,t^{-1}])$}, Algebr. Geom. Topol., to appear.

\bibitem{knudson.book}
K.~Knudson, \emph{Homology of linear groups}, Progress in mathematics Vol. 193; Basel, Boston, Berlin: Birkh$\ddot{a}$user, 2001.

\bibitem{knudson4}
K.~Knudson, \emph{Unstable homotopy invariance and the homology of $SL(2,\Z[t])$}, J. Pure Appl. Algebra \textbf{148} (2000), 255--266.

\bibitem{krammer}
D.~Krammer, \emph{Braid groups are linear}, Ann. of Math. \textbf{155} (2002), 131--156.

\bibitem{lee.szczarba}
R.~Lee and R.~Szczarba, \emph{On the homology and cohomology of congruence subgroups}, Invent. math. \textbf{33} (1976), 15-53.

\bibitem{lopez}
J.~Lopez, \emph{Lie algebras associated to congruence subgroups}, Ph.D. Thesis, University of Rochester, 2010.

\bibitem{lubotzky}
A.~Lubotzky, \emph{A group theoretic characterization of linear groups}, J. Algebra \textbf{113} (1988), 207--214.

\bibitem{serre}
J.P.~Serre, \emph{Lie algebras and Lie groups}, 1964 lectures given at Harvard University, Second Edition, New York: Springer-Verlag, 1992.

\bibitem{shimura}
G.~Shimura, \emph{Introduction to the arithmetic theory of automorphic functions}, Publications of the Mathematical Society of Japan, \textbf{11}, Princeton: Princeton University Press, 1971.

\bibitem{stallings}
J.~Stallings, \emph{On torsion-free groups with infinitely many ends}, Annals of Mathematics \textbf{88} (1968), 312-334.

\bibitem{swan}
R.~Swan, \emph{Groups of cohomological dimension one}, Journal of Algebra \textbf{12} (1969), 585-610.

\end{thebibliography}
\end{document}